\documentclass[twoside,reqno]{amsart}

\usepackage{amsmath}
\usepackage{amssymb,latexsym,amstext,amsthm}
\usepackage{verbatim} 
\usepackage{eucal} 

\newtheorem{thm}[equation]{Theorem}

\newtheorem{cor}[equation]{Corollary}
\newtheorem{lem}[equation]{Lemma}
\newtheorem{exa}[equation]{Example}

\newtheorem{rem}[equation]{Remark}

\def\aos{a\otimes s}
\def\das{\dd(\aa\otimes S)}

\def\zm{\bbbz_m}
\def\ib{\bar{i}}
\def\qb{\bar{q}}
\def\rb{\bar{r}}
\def\jb{\bar{j}}
\def\kb{\bar{k}}
\def\tb{\bar{t}}
\def\sb{\bar{s}}
\def\0b{\bar{0}}
\def\otimest{\overleftarrow{\otimes}}

\def\totimes{{\overrightarrow{\otimes}}}
\def\epi{\ep(\ib)}
\def\epj{\ep(\jb)}
\def\eps{\ep(\sb)}

\def\ept{\ep(\tb)}
\def\soplus{\boldsymbol{\;\;\oplus\;\;}}

\def\andd{\quad\hbox{and}\quad}

\def\a{\alpha}

\def\sub{\subseteq}

\def\1k{\frac{1}{k}}

\def\d{\delta}

\def\sg{\sigma}

\def\quadd{\quad\quad}

\def\ad{\hbox{ad}}

\def\bbbz{{\mathbb Z}}

\def\aa{\mathcal A}

\def\dd{\mathcal D}
\def\ep{\epsilon}

\def\span{\hbox{span}}

\def\bb{\mathcal B}

\def\pb{\bar{p}}

\newcommand\Aut{\operatorname{Aut}}
\newcommand\charr{\operatorname{char}}

\newcommand\End{\operatorname{End}}

\newcommand\Mult{\operatorname{Mult}}
\newcommand\Der{\operatorname{{\mathcal D}er}}
\newcommand\one{\operatorname{1}}

\begin{document}
\setcounter{page}{1}

\title{Derivations of tensor product of algebras}

\author{Saeid Azam}
\address[Saeid Azam]
{Department of Mathematics\\ University of Isfahan\\Isfahan, Iran,
P.O.Box 81745-163} \email{saeidazam@yahoo.com, azam@sci.ui.ac.ir}

\begin{abstract} We prove a theorem about the derivation algebra of the tensor product of two
algebras. As an application, we determine the derivation algebra
of the fixed point algebra of the tensor product of two algebras,
with respect to the tensor product of two finite order
automorphisms of the involved algebras. These results generalize
some well-know theorems in the literature.
\end{abstract}
\maketitle
\centerline{\it Dedicated to Professor Bruce Allison on the
occasion of his sixtieth birthday}

\medskip
\setcounter{section}{-1}
\section{Introduction}\label{introduction} In 1969, R.~E.~Block
\cite{B} showed that the algebra of derivations of the tensor
product of two algebras (satisfying certain finite dimensionality
conditions) can be expressed in terms of the algebra of
derivations and the centroid of each of the involved algebras. In
1986, G.~Benkart and R.~V.~Moody \cite{BM} used this (with a new
proof) to establish several interesting results about the
derivation algebra of the fixed points of the tensor product of
two algebras, with respect to the tensor product of two finite
order automorphisms of the involved algebras. They applied their
results to determine the algebra of derivations of several
important classes of infinite dimensional Lie algebras, including
twisted and untwisted affine Kac--Moody Lie algebras \cite{K},
Virasoro algebras, and some subclasses of extended affine Lie
algebras (for information about extended affine Lie algebras see
\cite{AABGP}, \cite{BGK} and \cite{N2}).

All algebras we consider will be over a field $k$. We denote by
$\dd(\aa)$ and $C(\aa)$ the algebra of derivations and centroid of
an algebra $\aa$, respectively. Let $\aa$ be a perfect algebra and
$S$ be a commutative associative unital algebra. It is proved in
\cite[Theorem~7.1]{B} and \cite[Theorem~1.1]{BM} that if $\aa$ is
finite dimensional then
$$
\dd(\aa\otimes S)=\dd(\aa)\otimes S\soplus
C(\aa)\otimes\dd(S).\eqno{(1)}
$$
Therefore any derivation $d\in\dd(\aa\otimes S)$ can be
represented as $$ d=\sum_{i\in I}d_i\otimes s_i+\sum_{j\in
J}\gamma_j\otimes d'_j,\eqno{(2)} $$ where $d_i$'s are in
$\dd(\aa)$, $d'_j$'s are in $\dd(S)$, $\{s_i\}_{i\in I}$ is a
basis of $S$ and $\{\gamma_j\}_{j\in J}$ is a basis of $C(\aa)$.
Moreover, the $d_i$'s and $d'_j$'s are zero except for finitely
many $i\in I$ and $j\in J$, respectively. It turns out that the
main reason for using finite dimensionality in proving (1) is to
show that the natural map $\psi:C(\aa)\otimes S\rightarrow
C(\aa\otimes S)$ is an isomorphism (see (\ref{eqcent1}) for the
definition of $\psi$). However, one can show that under certain
less restrictive conditions $\psi$ remains an isomorphism. Indeed
this holds if $\aa$ is finitely generated as an algebra over $k$
or as a module over its centroid, or if $\aa$ is unital (see Lemma
\ref{cent1}). When $\psi$ is an isomorphism we are able to prove a
generalization of (1) in the sense that each derivation $d$ of
$\aa\otimes S$ can be expressed in the form (2), where
$\{d_i\}_{i\in I}$ and $\{d'_j\}_{j\in J}$ are {\it summable}
families in $\dd(\aa)$ and $\dd(S)$, respectively. Here by a
summable family of endomorphisms $\{f_i\}_{i\in I}$ of a vector
space $V$, we mean that for each $v\in V$, $f_i(v)=0$, except for
finitely many $i\in I$ (see Theorem~\ref{thm1}).

Let $\sg_1$ and $\sg_2$ be period $m$ automorphisms of $\aa$ and
$S$, respectively. Then $\sg_1$, $\sg_2$ and $\sg_1\otimes\sg_2$
induce $\zm$-gradings on $\aa$, $S$, $\aa\otimes S$ and also on
the algebra of derivations and centroids of these algebras. Then
\cite[Theorem~1.3]{BM} states that if $\aa$ is finite dimensional,
$k$ contains all $m^{\hbox{th}}$-roots of unity ($p\nmid m$ if
$\charr(k)=p>0$) and the homogeneous subspace of $S$ of degree $1$
contains a unit, then the restriction map
$$\pi:\big(\dd(\aa\otimes S)\big)_{\0b}\rightarrow\dd\big((\aa\otimes
S)_{\0b}\big)$$ is an isomorphism. This has some very nice
applications. We have been able to generalize this theorem to the
extent that it holds for all algebras such that the map $\psi$ is
an isomorphism. Our approach to the proof is different from
\cite{BM} and it corrects an inaccuracy which occurs in the
surjectivity part of Benkart-Moody's proof (Remark \ref{rem2}).

Our interest in the algebra of derivations of tensor product of
algebras arise from the study of so called {\it iterated} and {\it
multi loop algebras} (see \cite{ABP}). These algebras cover some
very interesting classes of algebras including centerless affine
Lie algebras and in general almost all {\it centerless Lie tori}
\cite{ABFP} (see \cite{N1} for the definition of Lie tori). We are
considering the algebra of derivations of iterated and multi-loop
algebras in an ongoing project.

The author would like to thank Professor Bruce Allison for many
delightful discussions and suggestions.  He also would like to
thank Professor Stephen Berman for drawing his attention to the
derivation of iterated loop algebras. Finally, the author would
like to thank the people of the Department of Mathematical and
Statistical Sciences, University of Alberta, for their hospitality
during his visit.

\medskip
\section{\bf Centroids and pfgc algebras}\label{centroids}
\setcounter{equation}{0}

Throughout this work we fix a field $k$ and two algebras $\aa$ and
$S$ over $k$. All other algebras also will be over $k$. The {\it
multiplication algebra} of $\aa$, denoted $\Mult(\aa)$ is the
subalgebra if endomorphisms of $\aa$ over $k$ generated by the
identity element and left and right multiplication by elements of
$\aa$. The {\it centroid} of $\aa$ is by definition the set of
endomorphisms of $\aa$ which commute with all elements of
$\Mult(\aa)$. That is
$$
C(A)=\{\gamma\in\End(\aa)\mid
\gamma(xy)=\gamma(x)y=x\gamma(y)\hbox{ for all }x,y\in\aa\}.
$$
Clearly $C(\aa)$ is a unital subalgebra of $\End(\aa)$. Then $\aa$
can be considered as a left $C(\aa)$-module by $\gamma\cdot
a=\gamma(a)$, $\gamma\in C(\aa)$, $a\in\aa$.

An algebra $\aa$ is called {\it perfect} if $\aa\aa=\aa$. The
centroid of any perfect algebra is commutative \cite[Ch.~X, \S~1,
Lemma~1]{J}. Following \cite{ABP} we call $\aa$ a {\it pfgc
algebra} if
\begin{itemize}
\item [(i)] $\aa\not=0$, \item [(ii)] $\aa$ is perfect, \item
[(iii)] $\aa$ is finitely generated as a module over $C(\aa)$.
\end{itemize}

For $s\in S$ let $L_s$ denote the left multiplication by $s$, and
consider the linear map
\begin{equation}\label{eqcent1}
\begin{array}{c}
\psi:C(\aa)\otimes S\longrightarrow C(\aa\otimes S)\vspace{2mm}\\
\gamma\otimes s\longmapsto\gamma\otimes L_s.
\end{array}
\end{equation}
The following lemma indicates certain situations in which $\psi$
is an isomorphism of algebras. The proof of parts (i) and (iii)
can be found in \cite[\S 2]{ABP}, where they are considering
algebras over a ring. The proof of all four parts can be found in
the recent work \cite[Proposition~2.19]{BN}. Since this a central
tool for the rest of this work and since its proof is short and
simple in our setting, we provide the proof for the convenience of
reader.
\begin{lem}\label{cent1}
Let $S$ be a unital commutative associative algebra. Suppose one
of the following holds:

\begin{itemize}
\item [(i)] $\aa$ is perfect and $S$ is finite dimensional.

\item [(ii)] $\aa$ is perfect and finitely generated over $k$.

\item [(iii)] $\aa$ is a pfgc algebra.

\item [(iv)] $\aa$ is unital.
\end{itemize}
Then the map $\psi$ defined by (\ref{eqcent1}) is an isomorphism
of associative algebras.
\end{lem}

\begin{proof} It is clear that $\psi$ is a homomorphism. Since
$\aa$ is perfect, $C(\aa)$ is commutative. Fix a basis
$\{s_i\}_{i\in I}$ of $S$. Then any element $T$ of $C(\aa)\otimes
S$ can be written uniquely in the form $T=\sum_{i\in
I}\gamma_i\otimes s_i$, $\gamma_i\in C(\aa)$. If $\psi(T)=0$, then
$\sum_{i\in I}\gamma_i(a)\otimes L_{s_i}(1)=0$ for all $a\in\aa$.
Thus $\gamma_i=0$ for all $i\in I$ and so $\psi$ is 1-1.

To see $\psi$ is onto let $\Gamma\in C(\aa\otimes S)$. Then for
$a\in\aa$, $\Gamma(a\otimes 1)=\sum_{i\in I}\gamma_i(a)\otimes
s_i$, where $\gamma_i\in C(\aa)$ and $\gamma_i(a)=0$ except for a
finitely many $i\in I$. We show that under either of conditions
(i)-(iv) in the statement there is a finite subset $I_0$ of $I$
such that
\begin{equation}\label{1}
\sum_{i\in I}\gamma_i(a)\otimes s_i=\sum_{i\in
I_0}\gamma_i(a)\otimes s_i\quadd\hbox{for all }a\in\aa.
\end{equation}
To see this take $I_0=I$ if (i) holds. If $A$, either as an
algebra over $k$ or as a module over $C(\aa)$, is generated by
elements $a_1,\ldots,a_m$, take
$$
I_0=\bigcup_{j=1}^{m}\{i\in I\mid \gamma_i(a_j)\not=0\}.
$$
Finally if $\aa$ is unital, take $I_0=\{i\in
I\mid\gamma_i(1)\not=0\}$. In either of these cases $I_0$ is
finite and using the facts that $\gamma_i\in C(\aa)$ and $C(\aa)$
is commutative, it is easy to see that (\ref{1}) holds.

Now $T:=\sum_{i\in I_0}\gamma_i\otimes s_i$ is an element of
$C(\aa)\otimes S$ and for any $a,a'\in\aa$ and $s\in S$
\[
\big(\sum_{i\in I_0}\gamma_i\otimes s_i\big)(aa'\otimes
s)=\big(\sum_{i\in I_0} \gamma_i(a)\otimes s_i\big)(a'\otimes
s)=\Gamma(a\otimes 1)(a'\otimes s)=\Gamma(aa'\otimes s),
\]
but $\aa$ is perfect so $\psi(T)=\Gamma$ and we are done.
\qed\end{proof}

\medskip
\section{\bf Derivations}\label{derivations}
\setcounter{equation}{0}

In this section we generalize a result of \cite[Theorem~7.1]{B},
regarding the algebra of derivations of tensor product of two
algebras (see also \cite[Theorem~1.1]{BM}).

 Let $\aa$ be an algebra and $\bb$ be a
subalgebra of $\aa$. By a {\it derivation} from $\bb$ into $\aa$
we mean a $k$-linear map $\d:\bb\rightarrow\aa$ such that
$$
\d(bb')=\d(b)b'+b\d(b')\quad\hbox{for all }b,b'\in\bb.
$$
Denote the space of all such derivations by $\dd(\bb,\aa)$. This
space is usually denoted by $\Der(\bb,\aa)$, however we are using
the abbreviated notation $\dd(\bb,\aa)$ since it will appear
frequently. If $\aa=\bb$, we simply write it as $\dd(\aa)$.

Let $S$ be a commutative, associative unital algebra. Then
$\aa\otimes S$ can be considered as a $S$-bimodule by
$s'\cdot(a\otimes s)=(a\otimes s)\cdot s'=a\otimes ss'$, for
$a\in\aa$, $s,s'\in S$. We note that this action associates with
the product on $\aa\otimes S$. Let $\dd_S(\aa\otimes S)$ denote
the subalgebra of $\dd(\aa\otimes S)$ consisting of $S$-module
derivations. So if $d\in\dd(\aa\otimes S)$, then
$d\in\dd_S(\aa\otimes S)$ if and only if $$ d(a\otimes
ss')=s'd(a\otimes s)=d(a\otimes s)s'$$ for all $a\in\aa$, $s,s'\in
S$.

We note that the map
\begin{equation}\label{eqiso}
\begin{array}{c}
\tau:\dd_S(\aa\otimes S)\longrightarrow \dd(\aa\otimes 1,\aa\otimes S)\vspace{2mm}\\
d\longmapsto d_{|_{\aa\otimes 1}}
\end{array}
\end{equation}
is a vector space isomorphism. Therefore we can transfer the Lie
algebra structure on $\dd_{S}(\aa\otimes S)$ to $\dd(\aa\otimes
1,\aa\otimes S)$ by
\[ [d_1,d_2]=\tau([\tau^{-1}(d_1),\tau^{-1}(d_2)]),
\]$d_1,d_2\in\dd(\aa\otimes 1,\aa\otimes S).$
Then as Lie algebras
\begin{equation}\label{eqder1}
\dd_S(\aa\otimes S)\cong\dd(\aa\otimes 1,\aa\otimes S).
\end{equation}
Using this isomorphism we identify $\dd(\aa\otimes 1,\aa\otimes
S)$ as a subalgebra of $\dd(\aa\otimes S)$.

Finally we set
\[\dd_{\aa\otimes 1}(\aa\otimes S)=\{\d\in\dd(\aa\otimes S)
\mid\d(\aa\otimes 1)=0\}.\] Clearly $\dd_{\aa\otimes 1}(\aa\otimes
S)$ is a Lie subalgebras of $\dd(\aa\otimes S)$.

\begin{lem}\label{der1}
Let $\aa$ be an algebra and $S$ be a commutative, associative
unital algebra. Then
\[
\dd(\aa\otimes S)=\dd_S(\aa\otimes S)\soplus\dd_{\aa\otimes
1}(\aa\otimes S)
\]
where the sum on the right is the direct sum of vector spaces.
\end{lem}

\begin{proof} Let $\d\in\das$ and define $d\in\End(\aa\otimes S)$ by
$$
d(\aos)=\d(a\otimes 1)s=s\d(a\otimes 1),\quad a\in\aa,\;s\in S.
$$
To see that $d\in\dd_S(\aa\otimes S)$, let $a,a'\in\aa$ and
$s,s'\in S$. Then
\begin{eqnarray*}
d\big((\aos)(a'\otimes s')\big)&=&\d(aa'\otimes 1)ss'\\
&=&(\d(a\otimes 1)(a'\otimes 1))ss'+((a\otimes 1)\d(a'\otimes 1))ss'\\
&=&\d(a\otimes 1)s(a'\otimes 1)s'+(a\otimes 1)s\d(a'\otimes 1)s'\\
&=&d(a\otimes s)(a'\otimes s')+d(a\otimes s)(a'\otimes s'),
\end{eqnarray*}
and
\[d((a\otimes s)s')=d(a\otimes ss')=\d(a\otimes 1)ss'=d(a\otimes s)s'.
\]

Clearly $d=\d$ on $\aa\otimes 1$ and so $(\d-d)(\aa\otimes 1)=0$.
So $\d=d+(\d-d)$ where $d\in\dd_S(\aa\otimes S)$ and
$\d-d\in\dd_{\aa\otimes 1}(\aa\otimes S)$. It is now easy to see
that the sum is direct.\qed
\end{proof}

\medskip
Let $\one$ denote the identity operator on $\aa$ and consider
$\one\otimes S$ as a subalgebra of $C(\aa)\otimes S$.

\begin{lem}\label{der2}
Assume that $\aa$ is perfect and $S$ is commutative, associative
and unital. If the map $\psi$ defined by (\ref{eqcent1}) is an
isomorphism then as vector spaces
\[D_{\aa\otimes 1}(\aa\otimes S)\cong\dd\big(\one\otimes S,C(\aa)\otimes S\big).
\]
In particular, taking this as an identification, we may consider
$\dd\big(\one\otimes S,C(\aa)\otimes S\big)$ as a Lie subalgebra
of $\dd(\aa\otimes S)$.
\end{lem}

\begin{proof} Since $\aa$ is perfect, $C(\aa)\otimes S$ is a commutative associative
unital algebra. Using $\psi$ we define the linear map
\begin{equation}\label{eqiso2}
\Phi:\dd\big(\one\otimes S,C(\aa)\otimes
S\big)\longrightarrow\dd_{\aa\otimes 1}(\aa\otimes S)
\end{equation}
by
\[\Phi(d)(a\otimes s)=\psi\big(d(\one\otimes s)\big)(a\otimes 1),
\]
$d\in\dd\big(\one\otimes S,C(\aa)\otimes S\big)$, $a\in\aa$, $s\in
S$. Since $d(\one\otimes 1)=0$ we have $\Phi(d)(a\otimes 1)=0$. We
now show that $\Phi(d)$ is a derivation of $\one\otimes S$ into
$C(\aa)\otimes S$. So let $a,a'\in\aa$ and $s,s'\in S$. Since
$\psi\big(d(\one\otimes s)\big)$ and $\psi(d(\one\otimes s'))$ are
in $C(\aa\otimes S)$ we have
\begin{eqnarray*}
\psi d(\one\otimes s)(aa'\otimes s')&=&\big(\psi(d(\one\otimes s)(a\otimes 1)\big)(a'\otimes s')\\
\psi d(\one\otimes s')(aa'\otimes s)&=&(a\otimes
s)\big(\psi(d(\one\otimes s')(a'\otimes 1)\big).
\end{eqnarray*}
Therefore
\begin{eqnarray*}
\Phi(d)(aa'\otimes ss')&=&\psi\big(d(\one\otimes ss')\big)(aa'\otimes 1)\\
&=&[\psi d(\one\otimes s)(\one\otimes L_s')+
(\one\otimes L_s)\psi d(\one\otimes s')](aa'\otimes 1)\\
&=&\psi d(\one\otimes s)(aa'\otimes s')+
(\one\otimes L_s)(\psi d(\one\otimes s')(aa'\otimes 1))\\
&=&\big(\psi(d(\one\otimes s)(a\otimes 1)\big)(a'\otimes s')+
(a\otimes s)\big(\psi(d(\one\otimes s')(a'\otimes 1)\big)\\
&=&\Phi(d)(a\otimes s)(a'\otimes s')+(a\otimes s)\Phi(d)(a'\otimes
s').
\end{eqnarray*}

Next we show that $\Phi$ is 1-1 and onto. For this we define an
inverse map for $\Phi$ as follows. Set
\[\begin{array}{c}
\Phi':\dd_{\aa\otimes 1}(\aa\otimes S)\longrightarrow
\dd\big(\one\otimes S,C(\aa)\otimes S\big),\vspace{2mm}\\
d'\longmapsto (\psi^{-1}\circ\ad d'\circ\psi)_{|_{\one\otimes S}}.
\end{array}
\]
We note that $[\dd(\aa\otimes S),C(\aa\otimes S]\sub C(\aa\otimes
S)$ and so $\Phi'$ is well-defined. First we show that
$\Phi\circ\Phi'$ is the identity map on $\dd_{\aa\otimes
1}(\aa\otimes S)$. So let $d'\in\dd_{\aa\otimes 1}(\aa\otimes S)$,
$a\in\aa$ and $s\in S$. Since $d'(\aa\otimes 1)=0$ we have
$$\big(\ad d'\circ\psi(\one\otimes s)\big)(a\otimes 1)=
[d',\one\otimes L_s](a\otimes 1)=d'(a\otimes s).$$ Therefore
\begin{eqnarray*}
\big((\Phi\circ\Phi')\big)(d')(a\otimes s)&=& \Phi(\psi^{-1}\circ\ad d'\circ\psi)(\aos)\\
&=& \psi\big((\psi^{-1}\circ\ad d'\circ\psi)(\one\otimes s)\big)(a\otimes 1)\\
&=&\psi\big(\psi^{-1}(d'(a\otimes s)\big)=d'(\aos).
\end{eqnarray*}
Finally, we show that $\Phi'\circ\Phi$ is the identity map on
$\dd\big(\one\otimes S,C(\aa)\otimes S\big)$. Let
$d\in\dd\big(\one\otimes S,C(\aa)\otimes S\big)$, $a\in\aa$ and
$s,s'\in S$. Then
\begin{eqnarray*}
[\Phi(d),\one\otimes L_s](a\otimes s')&=& \psi\big(d(\one\otimes
ss')\big)(a\otimes 1)-(\one\otimes L_s)
\psi\big(d(\one\otimes s')\big)(a\otimes 1)\\
&=&\psi\big(d(\one\otimes s)(\one\otimes s')+(\one\otimes s)d(\one\otimes s')\big)(a\otimes 1)\\
&&\hphantom{\psi\big(d(\one\otimes s)(\one\otimes
s')}-(\one\otimes L_s)
\psi\big(d(\one\otimes s')\big)(a\otimes 1)\\
&=&\psi d(\one\otimes s)(a\otimes s').
\end{eqnarray*}
Thus
\begin{eqnarray*}
(\Phi'\circ\Phi)(d)(\one\otimes s)&=& \big(\psi^{-1}\circ\ad \Phi(d)\circ\psi)(\one\otimes s)\\
&=&\big( \psi^{-1}\ad \Phi(d)\big)(\one\otimes L_s)\\
&=&\psi^{-1}\big(\psi(d)(\one\otimes s)\big)=d(\one\otimes s).
\end{eqnarray*}
This completes the proof.\qed
\end{proof}

To state our next result we need to introduce the notion of a
summable family of endomorphisms on a vector space. A family
$\{f_i\}_{i\in I}$ of endomorphisms of a $k$-vector space $V$ is
called {\it summable} on $V$, if for each $v\in V$, $f_i(v)=0$
except for finitely many $i\in I$. If $\{f_i\}_{i\in I}$ is
summable, then we define $\sum_{i\in I}f_i\in\End(V)$ by
$$ \big(\sum_{i\in I}f_i\big)(v)=\sum_{i\in
I}f_i(v)\qquad\qquad(v\in V). $$

Now let $V$ and $W$ be two vector spaces over $k$ and let
$\hbox{E}_V$ and $\hbox{E}_W$ be subspaces of $\End(V)$ and
$\End(W)$ respectively.
If $\{f_i\}_{i\in I}\sub \hbox{E}_V$ is summable on $V$ and
$\{g_i\}_{i\in I}$ is any family in $\hbox{E}_W$, then
$\{f_i\otimes g_i\}_{i\in I}$ is summable on $V\otimes W$.
Therefore $\sum_{i\in I}f_i\otimes g_i\in\End(V\otimes W)$. We set
$$\hbox{E}_V\otimest \hbox{E}_W:=\{\sum_{i\in I}f_i\otimes g_i\mid\{f_i\}_{i\in
I}\sub \hbox{E}_V \hbox{ summable on }V,\;\{g_i\}_{i\in I}\sub
\hbox{E}_W\}.
$$
(The arrow points toward the subspace which the summable families
belong to.) Let $\{h_j\}_{j\in J}$ be a fixed basis of
$\hbox{E}_W$ and $\sum_{i\in I}f_i\otimes g_i\in
\hbox{E}_V\otimest \hbox{E}_W$. Let $g_i=\sum_{j\in J}\a^i_jh_j$,
$\a^i_j\in k$ and set $\tilde{f}_j=\sum_{i\in I}\a^i_jf_i$. Then
it is not difficult to show that $\{\tilde{f}_j\}_{j\in J}\sub
\hbox{E}_V$ is summable on $V$ and
$$
\sum_{i\in I}f_i\otimes g_i=\sum_{j\in J}\tilde{f}_j\otimes h_j.
$$
Therefore
$$
\hbox{E}_V\otimest \hbox{E}_W=\{\sum_{j\in J}f_j\otimes h_j\mid
\{f_j\}_{j\in J}\sub \hbox{E}_V\hbox{ summable on }V\}.
$$
Clearly $\hbox{E}_V\otimest \hbox{E}_W$ is a subspace of
$\End(V\otimes W)$. Similarly, if $\{f_i\}_{i\in I}$ is a basis of
$\hbox{E}_V$, we can define
$$
\hbox{E}_V\totimes \hbox{E}_W:=\{\sum_{i\in I}f_i\otimes g_i\mid
\{g_i\}_{i\in I}\sub \hbox{E}_W\hbox{ summable on }W\}.
$$

\medskip
We are now ready to state our next lemma. Note that when $S$ is
commutative and associative, we may identify it with a subspace of
$\End(S)$ through left (or right) multiplication.

\begin{lem}\label{tensor}
Let $\aa$ be perfect and $S$ be commutative associative and
unital. Then $\dd(\aa)\otimest S$ and $C(\aa)\totimes\dd(S)$ are
subalgebras of $\dd(\aa\otimes S)$. Moreover, if $\{s_i\}_{i\in
I}$ is a basis of $S$ and $\{\gamma_{j}\}_{j\in J}$ is a  basis of
$C(\aa)$ and $d\in\dd(\aa\otimes 1,\aa\otimes S)$,
$d'\in\dd(\one\otimes S,C(\aa)\otimes S)$, then for $a\in\aa$,
$s\in S$,
$$d(a\otimes 1)=\sum_{i\in I}d_i(a)\otimes
s_i\andd d'(\one\otimes s)=\sum_{j\in J}\gamma_j\otimes d'_j(s),
$$
where $\{d_i\}_{i\in I}$ and $\{d'_j\}_{j\in J}$ are summable
families in $\dd(\aa)$ and $\dd(S)$, respectively. In particular,
$$ \dd(\aa\otimes 1,\aa\otimes S )\cong\dd(\aa)\otimest S
\andd\dd\big(\one\otimes S,C(\aa)\otimes S\big)\cong
C(\aa)\totimes\dd(S),
$$
under the assignments
$$d\longmapsto\sum_{i\in I}d_i\otimes s_i\andd
d'\longmapsto\sum_{j\in J}\gamma_j\otimes d'_j,
$$
respectively.
\end{lem}

\begin{proof} We first show that $\dd(\aa)\otimest S$ is a subalgebra of
$\dd(\aa\otimes S)$. The result for $C(\aa)\totimes\dd(S)$ follows
by symmetry. Clearly $\dd(\aa)\otimest S$ is a subspace of
$\dd(\aa\otimes S)$. Now let $\{s_i\}_{i\in I}$ be a basis of $S$
and let $\{d_i\}_{i\in I}$ and $\{d'_i\}_{i\in I}$ be two summable
families in $\dd(\aa)$. Let $d=\sum_{i\in I}d_i\otimes s_i$ and
$d'=\sum_{i\in I}d'_i\otimes s_i$. Then for $a\in\aa$ and $s\in
S$, we have
$$
[d,d'](a\otimes s)=\sum_{i,j}[d_j,d'_i](a)\otimes s_is_js.
$$
Let $s_is_j=\sum_{t}\a^t_{i,j}s_t$. Since $\{d_i\}$ and $\{d'_i\}$
are summable, we have for each $t$ that
$\Delta_t=\sum_{i,j}\a^{t}_{i,j}[d_j,d'_i]$ is a well-defined
element of $\dd(\aa)$. Then
$$[d,d'](a\otimes s)=\sum_{t}\Delta_t(a)\otimes s_ts.$$
So we are done if we show that $\{\Delta_t\}_{t\in I}$ is
summable. Let $a\in\aa$ and set
$$
\begin{array}{l}
I_1=\{i\in I\mid d'_i(a)\not=0\},\\
I_2=\{i\in I\mid
d_i(a)\not=0\},\\
I_3=\{i\in I\mid d_i\big(d'_j(a)\big)\not=0,\;\;j\in I_1\cup
I_2\},\\I_4=\{i\in I\mid d'_i\big(d_j(a)\big)\not=0,\;\; j\in
I_1\cup I_2\},\\
I_0=\{t\in I\mid \a_{i,j}^t\not=0,\;i,j\in I_1\cup I_2\cup I_3\cup
I_4 \}.
\end{array}
$$
Now $I_0$ is finite and if $t\not\in I_0$, then for each $i,j$,
either $\a^t_{i,j}=0$ or $[d_j,d'_i](a)=0$. This shows that
$\Delta_t(a)=0$ if $t\not\in I_0$. Thus $\{\Delta_t\}_{t\in I}$ is
summable. This completes the proof of the first statement. The
proof of the second statement is straightforward.\qed
\end{proof}

\medskip
Recall that the {\it differential centroid} of an algebra $\aa$,
denoted by $dC(\aa)$, is by definition the centralizer of
$\dd(\aa)$ in $C(\aa)$, that is
$$dC(\aa)=\{\gamma\in C(\aa)\mid [\gamma,\dd(\aa)]=0\}.$$
We may consider $\aa$ also as a $dC(\aa)$-module.

\begin{lem}\label{der3}
Let $\aa$ be perfect and $S$ be commutative associative and
unital. Assume one of the following holds.
\begin{itemize}
\item [(a)] $S$ is finite dimensional, \item [(b)] $\aa$ is
finitely generated, \item [(c)] $\aa$ is a pfgc algebra and
$[\dd(\aa), C(\aa)]=0$, \item [(d)] $\aa$ is finitely generated as
a module over its differential centroid.
\end{itemize}
Then
$$\dd(\aa)\otimest S= \dd(\aa)\otimes S\andd
S\totimes\dd(\aa)=S\otimes\dd(\aa).
$$
\end{lem}

\begin{proof}
By symmetry, we only need to prove the first equality. We must
show that under either of conditions (a)-(d), any summable family
in $\dd(\aa)$ is finite. This is clear in the case (a). For the
cases (b)-(d) let $\{a_1,\ldots,a_m\}$ be a set of generators for
$\aa$ either as an algebra, as a module over $C(\aa)$ or as a
module over $dC(\aa)$, respectively. Let $\{d_i\}_{i\in I}$ be an
summable family in $\dd(\aa)$ and set
\[
I_0=\bigcup_{j=1}^{m}\{i\in I\mid d_i(a_j)\not=0\}.
\]
Then $I_0$ is finite. Let $i\in I\setminus I_0$. For the case (b)
it is clear that $d_i(\aa)=0$ and so $\{d_i\}_{i\in I}$ is finite.
For the case (c) we have
$$d_i(\aa)=\sum_{j=1}^{m}d_i\big(C(\aa)\cdot a_j\big)=\sum_{j=1}^{m}[d_i,C(\aa)](a_j)+\sum_{j=1}^{m}
C(\aa)d_i(a_j)=0.
$$
Argument for the case (d) is exactly the same as (c), replacing
the role of $C(\aa)$ by $dC(\aa)$.\qed
\end{proof}

\medskip
The following theorem summarizes the results of this section.

\begin{thm}\label{thm1} Let $\aa$ and $S$ be algebras such that $\aa$ is perfect
and $S$ is commutative associative and unital. Assume that the map
$\psi$ defined by (\ref{eqcent1}) is an isomorphism. Then
\begin{equation}\label{eqthm1}
\dd(\aa\otimes S)=\dd(\aa)\otimest S\soplus C(\aa)\totimes\dd(S).
\end{equation}
In particular, if either of the following is satisfied
\begin{itemize}
\item [(i)] $S$ is finite dimensional, \item [(ii)] $\aa$ is
finitely generated, \item [(iii)] $\aa$ is a pfgc algebra, \item
[(iv)] $\aa$ is unital,
\end{itemize}
then (\ref{eqthm1}) holds. Finally, if one of the following is
satisfied:
\begin{itemize}
\item [(i)] $\aa$ or $S$ is finite dimensional. \item [(ii)]
$C(\aa)$ is finite dimensional or $S$ is finitely generated over
$k$, and either of the following holds:
\begin{itemize}
\item [(a)] $\aa$ is finitely generated over $k$, \item [(b)]
$\aa$ is finitely generated as a $C(\aa)$-module and
$[\dd(\aa),C(\aa)]=0$, \item [(c)]$\aa$ is finitely generated as a
$dC(\aa)$-module,
\end{itemize}
\end{itemize}
then
\begin{equation}\label{eqthm2}
\dd(\aa\otimes S)=\dd(\aa)\otimes S\soplus C(\aa)\otimes \dd(S).
\end{equation}
\end{thm}

\begin{proof} The first statement follows from
(\ref{eqder1}) and Lemmas \ref{der1}, \ref{der2}, \ref{tensor}.
The second statement then is clear from Lemma \ref{cent1}.
Finally, the last statement follows from Lemma \ref{der3} and the
previous statements. \qed
\end{proof}

\medskip
\begin{rem}\label{rem1}
{\em Theorem \ref{thm1} is a generalization of \cite[Theorem
7.1]{B} and \cite[Theorem 1.1]{BM}. In fact when $\aa$ is finite
dimensional Theorem \ref{thm1} coincides with those of \cite{B}
and \cite{BM} (the work of \cite{B} contains more results when
$\aa$ is unital).}
\end{rem}

\medskip
\section{\bf Gradings induced by automorphisms}\label{gradations}
\setcounter{equation}{0} In this section we discuss the gradations
induced by finite order automorphisms. We assume that $k$ contains
a primitive $m^{\hbox{th}}$-root of unity. Starting from two
period $m$ automorphisms $\sg_1$ and $\sg_2$ of $\aa$ and $S$,
respectively, we consider the various $\zm$-gradations they induce
and investigate the relation between these gradations.

Let $\zm$ be the group of integers congruent to $m$ and let
$i:\bbbz\rightarrow \zm$ be the canonical map. Recall that a
$\zm$-{\it grading} of the algebra $\aa$ is an indexed family
$\{\aa_{\ib}\}_{\ib\in\zm}$ of subspaces of $\aa$ so that
$\aa=\oplus_{\ib\in\zm}\aa_{\ib}$ and
$\aa_{\ib}\aa_{\jb}\sub\aa_{\ib+\jb}$ for $\ib,\jb\in\zm$. Let
$\sg$ be an automorphism of $\aa$ of period $m$, and let
$\aa=\sum_{\ib\in\zm}\aa_{\ib}$ be the corresponding gradation on
$\aa$ where
$$
\aa_{\ib}=\{a\in\aa\mid \sg(a)=\omega^i a\},
$$
$\omega$ a primitive $m^{th}$-root of unity. The automorphism
$\sg$ extends to an automorphism $\sg^*$ of $\End(\aa)$ of period
$m$, by
$$
\sg^*(T)=\sg T\sg^{-1}.
$$
This in turn induces a $\zm$-grading
$\End(\aa)=\sum_{\ib\in\zm}\End(\aa)_{\ib}$, where
$$
\End(\aa)_{\ib}=\{T\in\End(\aa)\mid
T(\aa_{\jb})\sub\aa_{\ib+\jb}\hbox{ for all }\jb\in\zm\}.
$$
By restriction $\sg^*$ induces two automorphisms on $\dd(\aa)$ and
$C(\aa)$, of period $m$, and so we also have the $\zm$-gradings
$\dd(\aa)=\sum_{\ib\in\zm}\dd(\aa)_{\ib}$ and
$C(\aa)=\sum_{\ib\in\zm}C(\aa)_{\ib}$ where
$$
\dd(\aa)_{\ib}=\End(\aa)_{\ib}\cap\dd(\aa)\andd
C(\aa)_{\ib}=\End(\aa)_{\ib}\cap C(\aa).
$$

Next let $\sg_1$ and $\sg_2$ be automorphisms of period $m$ of
algebras $\aa$ and $S$, respectively. Then $\sg=\sg_1\otimes
\sg_2$ is a period $m$ automorphism of $\aa\otimes S$. So
$\aa\otimes S=\sum_{\ib\in\zm}(\aa\otimes S)_{\ib}$ where
$(\aa\otimes S)_{\ib}=\sum_{\jb\in\zm}\aa_{\ib-\jb}\otimes
S_{\jb}$. Also, by restriction, $\sg^*$ induces automorphisms of
$$\dd_S(\aa\otimes S)
\andd \dd_{\aa\otimes 1}(\aa\otimes S).
$$

\begin{lem}\label{der5} Assume that the field $k$ contains all $m^{th}$-roots
of unity for some integer $m$, $\aa$ is perfect, $S$ is
commutative associative and unital, and that the map $\psi$
defined by (\ref{eqcent1}) is an isomorphism. Let $\sg_1$ and
$\sg_2$ be period $m$ automorphisms of $\aa$ and $S$,
respectively. Then with respect to the $\zm$-gradings induced from
$\sg^*=\sg_1^*\otimes\sg_2^*$, $\sg_1^*$ and $\sg_2^*$ we have
$$\dd(\aa\otimes S)=\sum_{\jb\in\zm}\dd(\aa\otimes
S)_{\jb}=\sum_{\jb\in\zm}\big(\dd(\aa)\otimest S\big)_{\jb}\soplus
\big(C(\aa)\totimes\dd(S)\big)_{\jb},
$$
where
$$\big(\dd(\aa)\otimest S\big)_{\jb}=\sum_{\kb\in\zm}\dd(\aa)_{\kb}\otimest
S_{\jb-\kb}\andd
\big(C(\aa)\totimes\dd(S\big)_{\jb}=\sum_{\kb\in\zm}C(\aa)_{\kb}\totimes
\dd(S)_{\jb-\kb}.
$$
\end{lem}

\begin{proof} The first two equalities in the statement
follows from Lemma \ref{der1} and (\ref{eqthm1}) of Theorem
\ref{thm1}, where under the identification made in Section
\ref{derivations}, we have
\begin{equation}\label{eqder6}
\big(\dd(\aa)\otimest S\big)_{\jb}=\{d\in\dd(\aa)\otimest S\mid
d(\aa_{\ib}\otimes 1)\sub (\aa\otimes S)_{\ib+\jb}\hbox{ for all
}\ib\in\zm\}\vspace{2mm},
\end{equation}
and
$$ \big(C(\aa)\totimes\dd(S)\big)_{\jb}=\{d\in
C(\aa)\totimes\dd(S)\mid d(\one\otimes S_{\ib})\sub
\big(C(\aa)\otimes S\big)_{\ib+\jb}\hbox{ for all }\ib\in\zm\}.
$$

We continue the proof by proving the third equality, the fourth
equality follows by symmetry. From (\ref{eqder6}) it is clear that
the right hand side of this equality is a subset of the left hand
side. To see the reverse inclusion, let $d\in\big(\dd(\aa)\otimes
S\big)_{\jb}$ and let $\{s_i\}_{i\in I}$ be a basis of $S$
consisting of homogeneous elements with respect o the
$\zm$-grading on $S$. For $\kb\in\zm$ let $I_{\kb}=\{i\in I\mid
s_i\in S_{\kb}\}$. Now $d=\sum_{i\in I}d_i\otimes s_i$ where
$\{d_i\}_{i\in I}$ is an summable family in $\dd(\aa)$. For
$\kb\in\zm$ define
$$d^{\kb}_i=\left\{\begin{matrix}
d_i &\hbox{if }i\in I_{\kb}\\
0&\quad\hbox{otherwise}.
\end{matrix}\right.
$$
Then $\{d^{\kb}_i\}_{i\in I}$ is an summable family in $\dd(\aa)$.
Moreover, we have $d=\sum_{\kb\in\zm}d^{\kb}$ where
$d^{\kb}=\sum_{i\in I}d^{\kb}_i\otimes s_i.$ So it is enough to
show that for each $\kb\in\zm$,
$$d^{\kb}\in\sum_{\kb\in\zm}\dd(\aa)_{\kb}\otimest S_{\jb-\kb}.$$
Therefore without loss of generality we may assume that
$d=\sum_{i\in I}d^{\kb}_i\otimes s_i$ for some $\kb\in\zm$. This
means that $d^{\kb}_i\otimes s_i\in\dd(\aa)\otimes S_{\kb}$ for
all $i\in I$. Next, with respect to the $\zm$-grading on
$\dd(\aa)$ we have $d^{\kb}_i=\sum_{\tb\in\zm}d^{\kb,\tb}_i$ where
$d^{\kb,\tb}_i\in\dd(\aa)_{\tb}$. So $$
d=\sum_{\tb\in\zm}\sum_{i\in I}d^{\kb,\tb}_i\otimes s_i, $$ where
for each $\tb\in\zm$, the family $\{d^{\kb,\tb}_i\}_{i\in I}$ is
summable in $\dd(\aa)_{\tb}$ and $d^{\kb,\tb}_i\otimes
s_i\in\dd(\aa)_{\tb}\otimes S_{\kb}$ for all $i\in I$. So again
without loss of generality we may assume that $d=\sum_{i\in
I}d^{\kb,\tb}_i\otimes s_i$ where for each $i$,
$d^{\kb,\tb}_i\otimes s_i\in\dd(\aa)_{\tb}\otimes S_{\kb}$. Now
since $d\in\big(\dd(\aa)\otimes S\big)_{\jb}$ we have from
(\ref{eqder6}) that for each $\bar{l}\in\zm$,
$d(\aa_{\bar{l}}\otimes 1)\sub (\aa\otimes S)_{\bar{l}+\jb}$. But
for each $i$,
$$(d^{\kb,\tb}_i\otimes s_i)
(\aa_{\bar{l}}\otimes 1)\sub \aa_{\tb+\bar{l}}\otimes
S_{\kb}\sub(\aa\otimes S)_{\kb+\bar{l}+\tb},$$ and so
$\kb+\tb=\jb$. So for each $i$, $d^{\kb,\tb}_i\otimes
s_i\in\dd(\aa)_{\tb}\otimes S_{\jb-\tb}$ and we are done.\qed
\end{proof}

\bigskip
\section{\bf The interaction of fixed points and
derivations}\label{interaction} \setcounter{equation}{0} This
section contains the main result of this work (Theorem \ref{thm2})
which is a generalization of Theorem~1.3 of the very interesting
article \cite{BM}. The proof of Theorem \ref{thm2} in part
corrects the proof of \cite[Theorem~1.3]{BM} (see Remark
\ref{rem2}).

\begin{thm}\label{thm2} Let $\aa$ and $S$ be algebras over $k$ where $k$
contains all $m^{th}$-roots of unity for some integer $m$. If
$\charr(k)=p>0$ assume that $p \nmid m$.  Assume $\aa$  and $S$
satisfy the followings:
\begin{itemize}
\item [(i)] $\aa$ is perfect,

\item [(ii)] $S$ is commutative associative and unital,

\item [(iii)] $\sg_1\in\Aut(\aa)$, $\sg_2\in\Aut(S)$,
$\sg_1^m=\one$, and $\sg_2^m=\one$,

\item [(iv)] For some unit $q\in\zm$, there is a unit $u$ in
$S_q$.

\item [(v)] The map $\psi$ defined by (\ref{eqcent1}) is an
isomorphism,
\end{itemize}

If $(\aa\otimes S)_{\0b}$ denotes the fixed points of $\aa\otimes
S$ with respect to $\sg:=\sg_1\otimes\sg_2$ then the restriction
map
\begin{equation}\label{pi}
\begin{array}{c} \pi:(\dd(\aa\otimes
S)\big)_{\0b}\longrightarrow\dd\big((\aa\otimes
S)_{\0b}\big)\vspace{2mm}\\
\qquad D\longmapsto D_{|_{(\aa\otimes S)_{\0b}}}
\end{array}
\end{equation}
is an isomorphism. In particular,
\begin{equation}\label{thm2new}
\dd\big((\aa\otimes
S)_{\0b}\big)\cong\sum_{\ib\in\zm}\dd(\aa)_{\ib}\otimest
S_{-\ib}\soplus C(\aa)_{\ib}\totimes\dd(S)_{-\ib}.
\end{equation}
\end{thm}

Before starting the proof we make an important remark followed by
an example.

\begin{rem}\label{rem2}
{\em (i) When $\aa$ is finite dimensional, condition (v) of
Theorem \ref{thm2} is automatically satisfied (see Lemma
\ref{cent1}). In this case Theorem \ref{thm2} is identical to
\cite[Theorem 1.3]{BM}.

(ii) While checking the proof of \cite[Theorem 1.3]{BM}, we
realized an inaccuracy which occurs in the surjectivity part of
the proof. In fact the proof is based on the claim that the
restriction map (\ref{pi}) is an isomorphism. To show that $\pi$
is surjective, the authors consider $d\in\dd\big((\aa\otimes
S)_{\0b}\big)$ and extend it to an element $D\in\End(\aa\otimes
S)$ as follows: Let $\ib\in\zm$ and $0\leq s<m$. Since $\qb$ is
unit in $\zm$ there is a unique $0\leq r<m$ such that $\sb=\qb\rb$
in $\zm$. Then for $x_{\ib}\in\aa_{\ib}$ and $b_{-\ib+\sb}\in
S_{-\ib+\sb}$, define
\begin{equation}\label{eqBM}
D(x_{\ib}\otimes b_{-\ib+\sb})=(\one\otimes
L_{u^r})d(x_{\ib}\otimes u^{-r}b_{-\ib+\sb})=u^{r}d(x_{\ib}\otimes
u^{-r}b_{-\ib+\sb}).
\end{equation}
They claim that $D\in\big(\dd(\aa\otimes S)\big)_{\0b}$.
Unfortunately this is not true, we have provided a counterexample
in Example \ref{exaBM}. In Lemma \ref{thanksgod}, by assuming that
$d\in\dd\big((\aa\otimes S)_{\0b}\big)$ has an extension
$D\in\big(\dd(\aa\otimes S)\big)_{\0b}$, we extract the right
formula for $D$ in terms of $d$. Then it takes quite a bit of
non-straightforward work to show that this formula really provides
a derivation (see Lemmas \ref{this3}, \ref{thanksgod} and
\ref{surjective}).
Finally, we should
mention that our approach to the proof of \ref{thm2} is different
form \cite{BM}. }
\end{rem}

\medskip
\begin{exa}\label{exaBM}
{\em Let $\aa=k\one$ and let $S=k[z^{\pm}]$ be the algebra of
Laurent polynomials in variable $z$. Let $\sg_1=\one$ and
$\sg_2(z^{n})=w^nz^n$ where $w=e^{2\pi\sqrt{-1}/4}$. Then
$\aa=\aa_{\0b}$ and $S=\oplus_{\ib\in\bbbz_4}S_{\ib}$ where
$S_{\ib}=z^{i}k[z^{\pm 4}]$. Then $\one\otimes
zd/{dz}\in\big(\dd(\aa\otimes S)\big)_{\0b}$. So $d:=(\one\otimes
zd/{dz})_{|_{(\aa\otimes S)_{\0b}}}\in\dd\big((\aa\otimes
S)_{\0b}\big).$ Let $D\in\End(\aa\otimes S)$ be as in
(\ref{eqBM}), where $\qb=\bar{1}$ and $u=z\in\bbbz_{\bar{1}}$. Now
$1\otimes z^5\in\aa_{\0b}\otimes S_{\0b+\sb}$ where $\sb=\bar{1}$.
So $r=1$, is the unique integer with $0\leq r<4$ such that
$\rb\qb=\sb$. Then
$$D(\one\otimes z^5)=z^1d(1\otimes z^{-1}z^5)=4(1\otimes
z^5).$$ A similar computation shows that
$$
D(1\otimes z^3)=z^3d(1\otimes z^{-3}z^3)=0\andd D(1\otimes
z^2)=z^2d(1\otimes z^{-2}z^2)=0.$$ Thus $D$ is not a derivation.}
\end{exa}

\medskip
To proceed with the proof of the theorem we need a few lemmas.

For $\ib\in\zm$ let $\ep(\ib)$ be the unique preimage of $\ib$ in
$\{0,1,\ldots,m-1\}$, under the map
$\bar{\;}:\bbbz\rightarrow\zm$. Then
$$ \ep(\ib+\jb)=\left\{
\begin{array}{ll}
\ep(\ib)+\ep(\jb)&\hbox{if }\ep(\ib)+\ep(\jb)<m\\
\ep(\ib)+\ep(\jb)-m&\hbox{if }\ep(\ib)+\ep(\jb)\geq m.
\end{array}\right.
$$

As it is mentioned in Remark \ref{rem2}, the expression
(\ref{eqBM}) defined by \cite{BM} is not the right way of
extending an element $d\in\dd\big((\aa\otimes S)_{\0b}\big)$ to an
element $D\in\dd(\aa\otimes S)$. In the following lemma we analyze
what would be the right way of extending $d$. Before that we note
that if $\qb$ and $u$ are as in Theorem \ref{thm2}, and if $\qb_1$
is such that $\qb\qb_1=\bar{1}$ then $u\in S_{\qb}$ if and only if
$u':=u^{\ep(\qb_1)}\in S_{\bar{1}}$. Thus, condition (iv) of
Theorem \ref{thm2} is equivalent to the condition:
\begin{itemize}
\item[(iv)$'$] $S_{\bar{1}}$ contains a unit.
\end{itemize}
From now on and for the sake of simplicity {\it we work with
(iv)$'$ instead of (iv)}.

\begin{lem}\label{rightder} Under the conditions of Theorem \ref{thm2}
(with (iv)$'$ in place of (iv)) let $d\in\dd\big((\aa\otimes
S)_{\0b}\big)$ be extended to an element $D\in\big(\dd(\aa\otimes
S)\big)_{\0b}$. Let $a_{\ib}\in\aa_{\ib}$ and $b_{-\ib+\sb}\in
S_{-\ib+\sb}$. Then the followings hold:

({\em a}) For any integer $n$ such that $n^{-1}$ makes sense we
have
\begin{equation}\label{rightder1}
\begin{array}{l}
D(a_{\ib}\otimes b_{-\ib+ \sb})=  d(a_{\ib}\otimes u^{-\ep(\sb)}b_{-\ib+\sb})u^{\ep(\sb)}\\
\qquad\qquad\;+ \eps(mn)^{- 1}u^{\epi}[u^{-mn}d(a_{\ib}\otimes
u^{{-\ep(\ib)}+mn})-d(a_{\ib}\otimes u^{{-\ep(\ib)}})]b_{-\ib
+\sb}.
\end{array}
\end{equation}
\normalsize In particular if $\charr(k)=0$, this holds for any
nonzero integer $n$.

({\em b}) If $\charr(k)=p$, then
$$D(a_{\ib}\otimes
b_{-\ib+\sb})=d(a_{\ib}\otimes u^{-pr}b_{-\ib+\sb})u^{pr},$$ where
$r=\ep(\sb\pb^{-1}).$
\end{lem}

\begin{proof}
(a) Let $d$ and $D$ be as in the statement. By Lemma \ref{der5},
$D=D_1+D_2$ where
$$
D_1\in\sum_{\jb\in\zm}\dd(\aa)_{\jb}\otimest S_{-\jb}\andd
D_2\in\sum_{\jb\in\zm}C(\aa)_{\jb}\totimes \dd(S)_{-\jb}.
$$
Therefore we can write $D_1=\sum_{i\in I}d_i\otimes s_i$ and
$D_2=\sum_{i\in J}\gamma_i\otimes d'_{i}$, where for each $i$,
$d_i\otimes s_i\in\dd(\aa)_{\jb}\otimes S_{-\jb}$, for some
$\jb\in\zm$ and similarly $\gamma_i\otimes d'_i\in
C(\aa)_{\jb}\otimes D(S)_{-\jb}$ for some $\jb\in\zm$. Let
$x_{-\ib}=u^{-\ep(\sb)}b_{-\ib+\sb}$. Then $x_{-\ib}\in S_{-\ib}$
and
\begin{eqnarray*}
D(a_{\ib}\otimes b_{-\ib+\sb})&=&D_1(a_{\ib}\otimes
u^{\ep(\sb)}x_{-\ib})+
D_2(a_{\ib}\otimes u^{\ep(\sb)}x_{-\ib})\\
&=&D_1(a_{\ib}\otimes x_{-\ib})u^{\ep(\sb)}+D_2(a_{\ib}\otimes
u^{\ep(\sb)})x_{-\ib}+D_2(a_{\ib}\otimes x_{-\ib})u^{\ep(\sb)}\\
&=& D(a_{\ib}\otimes x_{-\ib})u^{\ep(\sb)}+D_2(a_{\ib}\otimes
u^{\ep(\sb)})x_{-\ib}.
\end{eqnarray*}
Since $D(a_{\ib}\otimes x_{-\ib})=d(a_{\ib}\otimes x_{-\ib})$, we
obtain
\begin{equation}\label{*}
D(a_{\ib}\otimes b_{-\ib+\sb})=d(a_{\ib}\otimes
x_{-\ib})u^{\ep(\sb)}+D_2(a_{\ib}\otimes u^{\ep(\sb)})x_{-\ib}.
\end{equation}

Now let $n$ be an integer such that $n^{-1}$ makes sense (if
$\charr(k)=0$, $n$ could be any nonzero integer, and if
$\charr(k)=p$, $n$ could be any integer with $(n,p)=1$), then for
any $t\in\bbbz$,
\begin{equation}\label{this1}
D_2(a_{\ib}\otimes u^t)= tu^{t-1}D_2(a_{\ib}\otimes u)
=tn^{-1}u^{t-n}D_2(a_{\ib}\otimes u^n).
\end{equation}
Also,
\begin{eqnarray*}
d(a_{\ib}\otimes u^{{-\ep(\ib)}}u^{nm})&=&D(a_{\ib}\otimes
u^{{-\ep(\ib)}+nm})\\
&=& D(a_{\ib}\otimes
u^{{-\ep(\ib)}})u^{nm}+D_2(a_{\ib}\otimes u^{nm})u^{{-\ep(\ib)}}\\
&=&d(a_{\ib}\otimes u^{{-\ep(\ib)}})u^{nm}+D_2(a_{\ib}\otimes
u^{nm})u^{{-\ep(\ib)}}.
\end{eqnarray*}
Thus
\begin{equation}\label{this2}
D_2(a_{\ib}\otimes u^{nm})=d(a_{\ib}\otimes
u^{{-\ep(\ib)}+nm})u^{\epi}-d(a_{\ib}\otimes
u^{{-\ep(\ib)}})u^{nm+\epi}.
\end{equation}
Then from (\ref{this1}) and (\ref{this2}) we have
\begin{eqnarray*}
&&D_2(a_{\ib}\otimes
u^{\ep(\sb)})\\
&&\qquad=\eps(nm)^{-1}u^{\eps-nm}D_2(a_{\ib}\otimes
u^{nm})\\
&&\qquad=\eps(nm)^{-1}u^{\eps-nm}\big(u^{\epi}d(a_{\ib}\otimes
u^{{-\ep(\ib)}+nm})-d(a_{\ib}\otimes u^{{-\ep(\ib)}})u^{nm+\epi}\big)\\
&&\qquad=\eps(nm)^{-1}\big(u^{\epi+\eps-nm}d(a_{\ib}\otimes
u^{{-\ep(\ib)}+nm})-u^{\eps+\epi}d(a_{\ib}\otimes
u^{{-\ep(\ib)}})\big).
\end{eqnarray*}
Replacing this in (\ref{*}) we obtain
\begin{eqnarray*}
&&D(a_{\ib}\otimes u^{\ep(\sb)}x_{-\ib})=d(a_{\ib}\otimes
x_{-\ib})u^{\ep(\sb)}\\
 &&\qquad+\eps(nm)^{-1}u^{\epi+\eps}[u^{-nm}d(a_{\ib}\otimes
u^{{-\ep(\ib)}+nm})-d(a_{\ib}\otimes u^{{-\ep(\ib)}})]x_{-\ib}.
\end{eqnarray*}
Replacing $x_{-\ib}$ with $u^{-\ep(\sb)}b_{-\ib+\sb}$ we obtain
\begin{eqnarray*}
&&D(a_{\ib}\otimes b_{-\ib+\sb})=d(a_{\ib}\otimes
u^{-\ep(\sb)}b_{-\ib+\sb})u^{\ep(\sb)}\\
&&\qquad\qquad+\eps(nm)^{-1}u^{\epi}[u^{-nm}d(a_{\ib}\otimes
u^{{-\ep(\ib)}+nm})-d(a_{\ib}\otimes
u^{{-\ep(\ib)}})]b_{-\ib+\sb},
\end{eqnarray*}
where $n$ is any integer such that $n^{-1}$ makes sense.

\medskip
(b) Let $r$ be as in the statement, then
\begin{eqnarray*}
D_2(a_{\ib}\otimes u^{-pr}b_{-\ib+\sb})&=&D_2(a_{\ib}\otimes
u^{-pr})b_{-\ib+\sb}+D_2(a_{\ib}\otimes b_{-\ib+\sb})u^{-pr}\\
&=& -pru^{-pr-1}D_2(a_{\ib}\otimes
u)b_{-\ib+\sb}+D_2(a_{\ib}\otimes
b_{-\ib+\sb})u^{-pr}\\
&=&D_2(a_{\ib}\otimes b_{-\ib+\sb})u^{-pr}.
\end{eqnarray*}
Thus
\begin{eqnarray*}
d(a_{\ib}\otimes u^{-pr}b_{-\ib+\sb})&=&D_1(a_{\ib}\otimes
b_{-\ib+\sb})u^{-pr}+D_2(a_{\ib}\otimes
b_{-\ib+\sb})u^{-pr}\\
&=&D(a_{\ib}\otimes b_{-\ib+\sb})u^{pr}.
\end{eqnarray*}
$\quad$\qed
\end{proof}

\medskip
\begin{cor}\label{bruce1-1}
Under the conditions of Theorem \ref{thm2} the map $\pi$ is
injective.
\end{cor}

\begin{proof} Let $D\in\big(\dd(\aa\otimes S)\big)_{\0b}$ and
$d:=\pi(D)=0$. By Lemma \ref{rightder}, $D=0$ is the unique
extension of $d$.\qed\end{proof}

\medskip
\begin{lem}\label{this3}
Let $\aa$ and $S$ be algebras which satisfy conditions (i)-(iv) of
Theorem \ref{thm2} (with (iv)$'$ in place of (iv)). Let
$d\in\dd\big((\aa\otimes S)_{\0b}\big)$, $i,n\in\bbbz$ and
$a\in\aa_{\ib}$. Then the following formulas hold:
\begin{equation}\label{formula1}
{u^{-nm}d(a\otimes u^{-i+nm})+u^{nm}d(a\otimes
u^{-i-nm})=2d(a\otimes u^{-i})},
\end{equation}
\begin{equation}\label{formula2}
u^{-nm}d(a\otimes u^{-i+nm})+nu^md(a\otimes
u^{-i-m})=(1+n)d(a\otimes u^{-i}),
\end{equation}
\begin{equation}\label{formula3}
u^{-nm}d(a\otimes u^{-i+nm})-nu^{-m}d(a\otimes
u^{-i+m})=(1-n)d(a\otimes u^{-i}).
\end{equation}
In particular, for $a_{\ib}\in\aa_{\ib}$ and
$a_{\jb}\in\aa_{\jb}$,
\begin{equation}\label{formula4}
\begin{array}{l}
u^{\ep(\ib+\jb)}[u^{-m}d(a_{\ib}a_{\jb}\otimes
u^{-\ep(\ib+\jb)+m})-d(a_{\ib}a_{\jb}\otimes
u^{-\ep(\ib+\jb)})]\vspace{2mm}\\
\qquad= u^{\ep(\ib)+\ep(\jb)}[u^{-m}d(a_{\ib}a_{\jb}\otimes
u^{-\ep(\ib)-\ep(\jb)+m})-d(a_{\ib}a_{\jb}\otimes
u^{-\ep(\ib)+\ep(\jb)})].
\end{array}
\end{equation}

\end{lem}

\begin{proof} We use induction on $n\geq 0$ to prove (\ref{formula1}).
Since $\aa$ is perfect, without loss of generality we may assume
that $a=xy$ where $x\in\aa_{\jb}$ and $y\in\aa_{\ib-\jb}$ for some
$\jb\in\zm$. Clearly formula (\ref{formula1}) holds for $n=0$. To
see it for $n=1$ note that
\begin{eqnarray*}
&&d(a\otimes u^{m-i})u^{-m}+d(a\otimes
u^{-m-i})u^{m}\\
&&\qquad\qquad=d(x\otimes u^{-j})(y\otimes
u^{-i+j})+(x\otimes u^{-j-m})d(y\otimes u^{-i+m+j})\\
&&\qquad\qquad\qquad+d(x\otimes u^{-j-m})(y\otimes
u^{-i+j+m})+(x\otimes u^{-j})d(y\otimes u^{-i+j})\\
&&\qquad\qquad=[d(x\otimes u^{-j})(y\otimes
u^{-i+j})+(x\otimes u^{-j})d(y\otimes u^{-i+j})]\\
&&\qquad\qquad\qquad+[(x\otimes u^{-j-m})d(y\otimes
u^{-i+m+j})+d(x\otimes u^{-j-m})(y\otimes u^{-i+j+m})]\\
&&\qquad\qquad=2d(a\otimes u^{-i}).
\end{eqnarray*}
Next assuming $n\geq 2$, we have (using induction hypothesis)
{\small
\begin{eqnarray*} &&d(a\otimes u^{nm-i})u^{-nm}+d(a\otimes
u^{-nm-i})u^{nm}\\
&&\qquad\quad=[d(x\otimes u^{m-j})(y\otimes
u^{-i+j+(n-1)m})+(x\otimes u^{-j+m})d(y\otimes u^{-i+j+(n-1)m})]u^{-nm}\\
&&\qquad\qquad\quad+[d(x\otimes u^{-j-m})(y\otimes
u^{-i+j-m(n-1)})+(x\otimes u^{-j-m})d(y\otimes u^{-i+j-(n-1)m})]u^{nm}\\
&&\qquad\quad=[d(x\otimes u^{-j+m})(y\otimes
u^{-i+j-m})+(x\otimes u^{-j-m(n-1)})d(y\otimes u^{-i+j+(n-1)m})]\\
&&\qquad\qquad\quad+[d(x\otimes u^{-j-m})(y\otimes
u^{-i+m+j})+(x\otimes u^{-j+(n-1)m})d(y\otimes u^{-i+j-(n-1)m})]\\
&&\qquad\quad=[d(x\otimes u^{-j+m})u^{-m}+d(x\otimes u^{-j-m})u^m](y\otimes u^{-i+j})\\
&&\qquad\qquad\quad +(x\otimes u^{-j})[u^{-(n-1)m}d(y\otimes
u^{-i+j+(n-1)m})+u^{m(n-1)}d(y\otimes u^{-i+j-m(n-1)}]\\
&&\qquad\quad=2d(x\otimes u^{-j})(y\otimes u^{-i+j})+2(x\otimes
u^{-j})d(y\otimes u^{-i+j})\\
&&\qquad\quad=2d(a\otimes u^{-i}).
\end{eqnarray*}}
\normalsize

Using induction, we first prove formula (\ref{formula2}) for
$n\geq 0$. It clearly holds for $n=0$ and it holds for $n=1$ by
(\ref{formula1}). So we may assume $n\geq 2$. Then using induction
hypothesis we have \small
\begin{eqnarray*}
&&d(a\otimes u^{nm-i})u^{-nm}+d(a\otimes
u^{-m-i})u^{m}\\
&&\qquad\qquad= d(x\otimes u^{-j+m})(y\otimes
u^{-i+j-m})+(x\otimes u^{-j-m(n-1)})d(y\otimes u^{j-i+m(n-1)})\\
&&\qquad\qquad\qquad+d(x\otimes u^{-j-m})(y\otimes
u^{-i+j+m})+(x\otimes u^{-j})d(y\otimes u^{-i+j})\\
&&\qquad\qquad=u^{-m}[d(x\otimes u^{-j+m})(y\otimes
u^{-i+j})+(x\otimes u^{-j+m})d(y\otimes u^{-i+j})]\\
&&\qquad\qquad\qquad+u^{-m(n-2)}[(x\otimes u^{-j-m})d(y\otimes
u^{j-i+m(n-1)})\\
&&\hspace{6cm}+d(x\otimes u^{-j-m})(y\otimes u^{-i+j+(n-1)m})]\\
&&\qquad\qquad=u^{-m}d(a\otimes
u^{m-i})+u^{-m(n-1)}d(a\otimes u^{-i+(n-2)m})\\
&&\qquad\qquad =u^{-m}d(a\otimes u^{m-i})+(n-1)d(a\otimes
u^{-i})-(n-2)u^md(a\otimes u^{-i-m})\\
&&\qquad\qquad=(n-1)d(a\otimes u^{-i})\\
&&\qquad\qquad\qquad+[u^{-m}d(a\otimes u^{m-i})+u^md(a\otimes
u^{-i-m})]-(n-1)u^md(a\otimes
u^{-i-m})\\
&&\qquad\qquad=(n-1)d(a\otimes u^{-i})+2d(a\otimes
u^i)-(n-1)u^md(a\otimes u^{-i-m})\\
&&\qquad\qquad=(n+1)d(a\otimes u^{-i})-(n-1)u^md(a\otimes
u^{-i-m}).
\end{eqnarray*}
\normalsize Thus (\ref{formula2}) holds for $n\geq 0$. To get
(\ref{formula2}) for $n\leq 0$, subtract (\ref{formula2}) for
$n\geq 0$ from (\ref{formula1}).

Next we prove (\ref{formula3}). From (\ref{formula1}) we have
$$nu^{-m}d(a\otimes u^{-i+m})+nu^{m}d(a\otimes
u^{-i-m})=2nd(a\otimes u^{-i}).$$ Now subtracting this from
(\ref{formula2}) we get (\ref{formula3}). Finally we show
(\ref{formula4}). If $\ep(\ib+\jb)=\epi+\epj$, there is nothing to
be proved. If $\ep(\ib+\jb)=\epi+\epj-m$, replace $a$ with
$a_{\ib}a_{\jb}$ and $i$ with $\ep(\ib)+\ep(\jb)$ in
(\ref{formula3}) for $n=2$ to get (\ref{formula4}).\qed
\end{proof}

\medskip
\begin{lem}\label{thanksgod}
Under the conditions (i)-(iv) of Theorem \ref{thm2} (with (iv)$'$
in place of (iv)), let $d\in\dd\big((\aa\otimes S)_{\0b}\big)$,
$i,j,s,t\in\bbbz$, $a_{\ib}\in\aa_{\ib}$, $a_{\jb}\in\aa_{\jb}$,
$b_{-\ib+\sb}\in S_{-\ib+\sb}$ and $b_{-\jb+\tb}\in S_{-\jb+\tb}$.
Then
\begin{equation}\label{*I}
\begin{array}{l}
[u^{-m+s}d(a_{\ib}\otimes u^{-s+m}b_{-\ib+
\sb})-u^sd(a_{\ib}\otimes
u^{-s}b_{-\ib+\sb})](a_{\jb}\otimes b_{-\jb+\tb})\\
\qquad\qquad= (a_{\ib}\otimes b_{-\ib+
\sb})[u^{-m+t}d(a_{\jb}\otimes u^{-t+m}b_{-\jb+\tb})-u^t
d(a_{\jb}\otimes u^{-t}b_{-\jb+\tb})].
\end{array}
\end{equation} In particular,
\begin{equation}\label{*II}
\begin{array}{l} (a_{\ib}\otimes u^{-i})[u^{-m+j}d(a_{\jb}\otimes
u^{-j+m})-d(a_{\jb}\otimes
u^{-j})u^j]u^{-j}\\
\qquad\qquad= [u^{-m+i}d(a_{\ib}\otimes u^{-i+m})-d(a_{\ib}\otimes
u^{-i})u^i]u^{-i}(a_{\jb}\otimes u^{-j}),
\end{array}
\end{equation}
and
\begin{equation}\begin{array}{l}\label{*III}
(a_{\ib}\otimes 1)[u^{-m+t}d(a_{\jb}\otimes u^{-t+m}b_{-\jb+\tb}
)-d(a_{\jb}\otimes
u^{-t}b_{-\jb+\tb})u^t]\\
\qquad\qquad = [u^{-m+i}d(a_{\ib}\otimes
u^{-i+m})-d(a_{\ib}\otimes u^{-i})u^i](a_{\jb}\otimes
b_{-\jb+\tb}).
\end{array}
\end{equation}
\end{lem}

\begin{proof}
For the sake of simplicity, we put $b_1=b_{-\ib+\sb}$ and
$b_2=b_{-\jb+\tb}$. Using the fact that $d$ is a derivation, we
compute $$M:=d(a_{\ib}a_{\jb}\otimes
u^{-s-t+m}b_1b_2)u^{s+t-m}-d(a_{\ib}a_{\jb}\otimes
u^{-s-t}b_1b_2)u^{s+t}$$ in the following two ways:
\begin{eqnarray*}
M&=& u^{-m+s+t}[d(a_{\ib}\otimes u^{-s+m}b_1)(a_{\jb}\otimes
u^{-t}b_2)
+(a_{\ib}\otimes u^{-s+m}b_1)d(a_{\jb}\otimes u^{-t}b_2)]\\
&&- u^{s+t}[d(a_{\ib}\otimes u^{-s}b_1)(a_{\jb}\otimes u^{-t}b_2)
+(a_{\ib}\otimes u^{-s}b_1)d(a_{\jb}\otimes u^{-t}b_2)]\\
&=& [d(a_{\ib}\otimes u^{-s+m}b_1)u^{s-m}(a_{\jb}\otimes b_2)
+(a_{\ib}\otimes b_1)d(a_{\jb}\otimes u^{-t}b_2)u^t]\\
&&-[d(a_{\ib}\otimes u^{-s}b_1)u^s(a_{\jb}\otimes b_2)
+(a_{\ib}\otimes b_1)d(a_{\jb}\otimes u^{-t}b_2)u^t]\\
&=& [u^{-m+s}d(a_{\ib}\otimes
u^{-s+m}b_{-\ib+\sb})-u^sd(a_{\ib}\otimes
u^{-s}b_{-\ib+\sb})](a_{\jb}\otimes b_{-\jb+\tb})
\end{eqnarray*}
and similarly
\begin{eqnarray*}
M&=&u^{-m+s+t}[d(a_{\ib}\otimes u^{-s}b_1)(a_{\jb}\otimes
u^{m-t}b_2)
+(a_{\ib}\otimes u^{-s}b_1)d(a_{\jb}\otimes u^{-t+m}b_2)]\\
&&- u^{s+t}[d(a_{\ib}\otimes u^{-s}b_1)(a_{\jb}\otimes u^{-t}b_2)
+(a_{\ib}\otimes u^{-s}b_1)d(a_{\jb}\otimes u^{-t}b_2)]\\
&=& (a_{\ib}\otimes b_{-\ib+\sb})[u^{-m+t}d(a_{\jb}\otimes
u^{-t+m}b_{-\jb+\tb})-u^t d(a_{\jb}\otimes u^{-t}b_{-\jb+\tb})].
\end{eqnarray*}
Now comparing the result of the above computations for $M$ we get
(\ref{*I}).

Next substitute $s=0$, $b_{-\ib+\sb}=u^{-i}$, $t=0$ and
$b_{-\jb+\tb}=u^{-j}$ in (\ref{*I}) to get (\ref{*II}). To get
(\ref{*III}) apply $s=i$ and $b_{-\ib+\ib}=1$ to (\ref{*I}).
\qed\end{proof}

\medskip
\begin{lem}\label{surjective}
Under the conditions of Theorem \ref{thm2} (with (iv)$'$ in place
of (iv)), the map $\pi$ defined by (\ref{pi}) is surjective.
\end{lem}

\begin{proof}
Let $d\in\dd\big((\aa\otimes S)_{\0b}\big)$. Define
$D\in\End(\aa\otimes S)$ by
\begin{eqnarray*}
&&D(a_{\ib}\otimes b_{-\ib+\sb})=d(a_{\ib}\otimes
u^{-\eps}b_{-\ib+\sb})u^{\eps}\\
&&\qquad\qquad\qquad+ \eps m^{- 1}u^{\epi}[u^{-m}d(a_{\ib}\otimes
u^{-\epi+m})-d(a_{\ib}\otimes u^{-\epi})]b_{-\ib+\sb},
\end{eqnarray*}
where $a_{\ib}\in\aa_{\ib}$ and $b_{-\ib+\sb}\in S_{-\ib+\sb}$. We
are done if we show that $D$ is a derivation of $\aa\otimes S$.
Let $a_{\ib}\in\aa_{\ib}$, $a_{\jb}\in\aa_{\jb}$,
$b_1=b_{-\ib+\sb}\in S_{-\ib+\sb}$ and $b_2=b_{-\jb+\tb}\in
S_{-\jb+\tb}$. Then by the definition of $D$ we have
\begin{eqnarray*}
A:&=&D(a_{\ib}a_{\jb}\otimes b_1b_2)=d(a_{\ib}a_{\jb}\otimes u^{-\ep(\sb+\tb)}b_1b_2)u^{\ep(\sb+\tb)}\\
&&+
\ep(\sb+\tb)m^{-1}u^{\ep(\ib+\jb)}[u^{-m}d(a_{\ib}a_{\jb}\otimes
u^{-\ep(\ib+\jb)+m})-d(a_{\ib}a_{\jb}\otimes
u^{-\ep(\ib+\jb)})]b_1b_2.
\end{eqnarray*}
By (\ref{formula4}) we can change $\ep(\ib+\jb)$ to
$\ep(\ib)+\ep(\jb)$ in the above expression, then we obtain
\begin{equation}\label{I}
\begin{array}{l}
A= d(a_{\ib}a_{\jb}\otimes
u^{-\ep(\sb+\tb)}b_1b_2)u^{\ep(\sb+\tb)}\\
\qquad\qquad\quad+\ep(\sb+\tb)m^{-1}u^{\epi+\epj}[u^{-m}d(a_{\ib}a_{\jb}\otimes
u^{-\epi- \epj+m})\\
\qquad\qquad\qquad\qquad\qquad\qquad\qquad\qquad\qquad\qquad-d(a_{\ib}a_{\jb}\otimes
u^{-\epi-\epj})]b_1b_2.
\end{array}
\end{equation}
We divide the argument to the cases $\ep(\sb+\tb)=\eps+\ept$ and
$\ep(\sb+\tb)=\eps+\ept-m$.

Case \underline{$\ep(\sb+\tb)=\ep(\sb)+\ep(\tb)$}: We have from
(\ref{I})
\begin{eqnarray*}
A&=&[d(a_{\ib}\otimes u^{-\eps}b_1)(a_{\jb}\otimes
u^{-\ept}b_2)+(a_{\ib}\otimes
u^{-\eps}b_1)d(a_{\jb}\otimes u^{-\ept}b_2)]u^{\eps+\ept}\\
&&+(\eps+\ept)m^{-1}u^{\epi+\epj} [u^{-m}d(a_{\ib}\otimes
u^{-\epi})(a_{\jb}\otimes
u^{-\epj+m})\\
&&\qquad\qquad\qquad\qquad\qquad\qquad\qquad+u^{-m}(a_{\ib}\otimes
u^{-\epi})d(a_{\jb}\otimes u^{-\epj+m})]b_1b_2\\
&&-(\eps+\ept)m^{-1}u^{\epi+\epj} [d(a_{\ib}\otimes
u^{-\epi})(a_{\jb}\otimes
u^{-j})\\
&&\qquad\qquad\qquad\qquad\qquad\qquad\qquad+(a_{\ib}\otimes u^{-\epi})d(a_{\jb}\otimes u^{-\epj})]b_1b_2\\
&=&d(a_{\ib}\otimes u^{-\eps}b_1)u^{\eps}(a_{\jb}\otimes b_2)+
(a_{\ib}\otimes
b_1)d(a_{\jb}\otimes u^{-\ept}b_2)u^{\ept}\\
&&+\eps m^{-1}(a_{\ib}\otimes b_1)[d(a_{\jb}\otimes
u^{-\epj+m})u^{\epj-m}-
d(a_{\jb}\otimes u^{-\epj})u^{\epj}]b_2\\
&&+ \ept m^{-1}(a_{\ib}\otimes b_1)[d(a_{\jb}\otimes
u^{-\epj+m})u^{\epj-m}- d(a_{\jb}\otimes u^{-\epj})u^{\epj}]b_2.
\end{eqnarray*}

On the other hand
\begin{eqnarray*}
B:&=&D(a_{\ib}\otimes b_1)(a_{\jb}\otimes b_2)+(a_{\ib}\otimes
b_1)D(a_{\jb}\otimes
b_2)\\
&=&\big(d(a_{\ib}\otimes u^{-\eps}b_1)u^{\eps}+\eps
m^{-1}u^{\epi}[u^{-m}d(a_{\ib}\otimes
u^{-\epi+m})\\
&&\qquad\qquad\qquad\qquad\qquad\qquad\qquad\qquad\qquad\qquad
-d(a_{\ib}\otimes u^{-\epi})]b_1\big)(a_{\jb}\otimes b_2)\\
&&+ (a_{\ib}\otimes b_1) \big(d(a_{\jb}\otimes
u^{-\ept}b_2)u^{\ept}+\ept m^{-1}u^{\epj}[u^{-m}d(a_{\jb}\otimes
u^{-\epj+m})\\
&&\qquad\qquad\qquad\qquad\qquad\qquad\qquad\qquad\qquad\qquad\qquad\qquad-d(a_{\jb}\otimes
u^{-\epj})]b_2\big)\\
&=&d(a_{\ib}\otimes u^{-\eps}b_1)u^{\eps}(a_{\jb}\otimes b_2)+
(a_{\ib}\otimes
b_1)d(a_{\jb}\otimes u^{-\ept}b_2)u^{\ept}\\
&&+\eps m^{-1}(a_{\ib}\otimes b_1)[u^{\epi-m}d(a_{\ib}\otimes
u^{-\epi+m})-d(a_{\ib}\otimes u^{-\epi})u^{\epi}]b_1(a_{\jb}\otimes b_2)\\
&&+\ept m^{-1}(a_{\ib}\otimes b_1)[u^{-m+\epj}d(a_{\jb}\otimes
u^{-\epj+m})-d(a_{\jb}\otimes u^{-\epj})u^{\epj}]b_2.
\end{eqnarray*}
We must show $A=B$. But from the above computations we see that
$A=B$ if and only if \begin{eqnarray*} &&\eps
m^{-1}(a_{\ib}\otimes b_1)[d(a_{\jb}\otimes
u^{-\epj+m})u^{\epj-m}-d(a_{\jb}\otimes u^{-\epj})u^{\epj}]b_2\\
&&\qquad=\eps m^{-1}[d(a_{\ib}\otimes
u^{-\epi+m})u^{\epi-m}-d(a_{\ib}\otimes
u^{-\epi})u^{\epi}]b_1(a_{\jb}\otimes b_2).
\end{eqnarray*}
But this holds if and only if
\begin{eqnarray*}
&&\eps m^{-1}(a_{\ib}\otimes u^{-\epi})[d(a_{\jb}\otimes
u^{-\epj+m})u^{\epj-m}\\
&&\qquad\qquad\qquad\qquad\qquad\qquad\qquad\qquad
-d(a_{\jb}\otimes u^{-\epj})u^{\epj}]u^{-\epj}b_1b_2\\
&&\qquad=\eps m^{-1}[d(a_{\ib}\otimes
u^{-\epi+m})u^{\epi-m}\\
&&\qquad\qquad\qquad\qquad\qquad\qquad-d(a_{\ib}\otimes
u^{-\epi})u^{\epi}](a_{\jb}\otimes u^{-\epj})u^{-\epi}b_1b_2.
\end{eqnarray*}
Now this holds by (\ref{*II}) of Lemma \ref{thanksgod}.

\medskip
Case \underline{$\ep(\sb+\tb)=\eps+\ept-m$}: In this case we have
from (\ref{I}) that
\begin{eqnarray*}
A&=&d(a_{\ib}a_{\jb}\otimes u^{-\eps-\ept+m}b_1b_2)u^{\eps+\ept-m}\\
&& \qquad\qquad
+(\eps+\ept-m)m^{-1}u^{\epi+\epj}[u^{-m}d(a_{\ib}a_{\jb}\otimes
u^{-\epi-\epj+m})\\
&&\qquad\qquad\qquad\qquad\qquad\qquad\qquad\qquad\qquad\qquad\qquad-d(a_{\ib}a_{\jb}\otimes
u^{-\epi-\epj})]b_1b_2
\end{eqnarray*}
Adding and subtracting the term $d(a_{\ib}a_{\jb}\otimes
u^{-\eps-\ept}b_1b_2)u^{\eps+\ept}$ to $A$, and applying the
result for the case $\ep(\sb+\tb)=\eps+\ept$, we obtain
$$
A=D(a_{\ib}\otimes b_1)(a_{\jb}\otimes b_2)+(a_{\ib}\otimes
b_1)D(a_{\jb}\otimes b_2)+M-N
$$
where
\begin{eqnarray*}
M=d(a_{\ib}a_{\jb}\otimes
u^{-\eps-\ept+m}b_1b_2)u^{\eps+\ept-m}-d(a_{\ib}a_{\jb}\otimes
u^{-\eps-\ept}b_1b_2)u^{\eps+\ept}
\end{eqnarray*}
and
$$N=u^{\epi+\epj}[u^{-m}d(a_{\ib}a_{\jb}\otimes
u^{-\epi-\epj+m})-d(a_{\ib}a_{\jb}\otimes
u^{-\epi-\epj})]b_1b_2.$$ So we are done if we show that $M-N=0$.
We have seen in the proof of Lemma \ref{thanksgod} that
\begin{eqnarray*}
M=(a_{\ib}\otimes b_1)[d(a_{\jb}\otimes
u^{-\ept+m}b_2)u^{\ept-m}-d(a_{\jb}\otimes u^{-\ept}b_2)u^{\ept}].
\end{eqnarray*}
Also we have
\begin{eqnarray*}
N&=&d(a_{\ib}\otimes u^{-\epi+m})u^{\epi-m}b_1(a_{\jb}\otimes
b_2)+(a_{\ib}\otimes
b_1)d(a_{\jb}\otimes u^{-\epj})u^{\epj}b_2\\
&&-d(a_{\ib}\otimes u^{-\epi})u^{\epi}b_1(a_{\jb}\otimes
b_2)-(a_{\ib}\otimes
b_1)d(a_{\jb}\otimes u^{-\epj})u^{\epj}b_2\\
&&=[d(a_{\ib}\otimes u^{-\epi+m})u^{\epi-m}-d(a_{\ib}\otimes
u^{-\epi})u^{\epi}]b_1(a_{\jb}\otimes b_2).
\end{eqnarray*}
Now by acting $b_1$ to the both sides of equality (\ref{*III}) of
Lemma \ref{thanksgod} we see that $M=N$. This completes the proof.
\qed
\end{proof}

\medskip
\noindent{\bf Proof of Theorem \ref{thm2}.} By Corollary
\ref{bruce1-1} and Lemma \ref{surjective} the map $\pi$ is an
isomorphism. The last part of the statement follows from Lemma
\ref{der5}. This completes the proof of theorem.\qed

\begin{cor}\label{cormain}
 Let $\aa$ and $S$ be algebras over $k$ where $k$
contains all $m^{th}$-roots of unity for some integer $m$. If
$\charr(k)=p>0$ assume that $p \nmid m$.  Assume $\aa$  and $S$
satisfy the followings:
\begin{itemize}
\item [(i)] $\aa$ is perfect and finitely generated,

\item [(ii)] $S$ is commutative associative unital and finitely
generated,

\item [(iii)] $\sg_1\in\Aut(\aa)$, $\sg_2\in\Aut(S)$,
$\sg_1^m=\one$, and $\sg_2^m=\one$,

\item [(iv)] For some unit $q\in\zm$, there is a unit $u$ in
$S_q$.

\end{itemize}
Then
$$\dd\big((\aa\otimes
S)_{\0b}\big)\cong\sum_{\ib\in\zm}\dd(\aa)_{\ib}\otimes
S_{-\ib}\soplus C(\aa)_{\ib}\otimes\dd(S)_{-\ib}.
$$
\end{cor}

\begin{proof}
By Lemma \ref{cent1}, condition (v) of Theorem \ref{thm2} is
satisfied and so (\ref{thm2new}) holds. Now the result follows
from Theorem \ref{thm1}.\qed
\end{proof}

\begin{rem}\label{remnew}{
(i) Corollary \ref{cormain} holds (with the same proof) replacing
condition (i) with one of the followings:
\begin{itemize}
\item [-] $\aa$ is a pfgc algebra with $[C(\aa),\dd(\aa)]=0$.
\item [-] $\aa$ is perfect and finitely generated as a module over
$dC(\aa)$.
\end{itemize}

(ii) Condition (iv) of Theorem \ref{thm2} (or Corollary
\ref{cormain}) can never be satisfied if $\sg_2=\hbox{id}$.
However if $\sg_1=\hbox{id}$ and $\sg_2=\hbox{id}$, this theorem
(or corollary) holds by Theorem \ref{thm1}. }
\end{rem}

\medskip
Under the conditions of Theorem \ref{thm2}, let
$\varphi:=\pi^{-1}$. Then Lemma \ref{rightder} gives the exact
formula for $\varphi$, namely for $d\in\dd\big((\aa\otimes
S)_{\0b}\big)$, $a_{\ib}\in\aa_{\ib}$ and $b_{-\ib+\sb}\in
S_{-\ib+\sb}$,
\begin{equation}\label{inversemap}
\begin{array}{l}
\varphi(d)(a_{\ib}\otimes b_{-\ib+ \sb})=  d(a_{\ib}\otimes u^{-\ep(\sb)}b_{-\ib+\sb})u^{\ep(\sb)}\\
\qquad\qquad\;+ \eps(m)^{- 1}u^{\epi}[u^{-m}d(a_{\ib}\otimes
u^{{-\ep(\ib)}+m})-d(a_{\ib}\otimes u^{{-\ep(\ib)}})]b_{-\ib
+\sb}.
\end{array}
\end{equation}

In Example \ref{exaBM} we saw that the map suggested in
\cite[Theorem~1.7]{BM} does not provide an inverse map for $\pi$.
In the following example we see how our inverse map $\varphi$
works for that particular example. We also explain how it works
for the most studied example, namely twisted affine Lie algebras.

\begin{exa}\label{last exa}
{\em (i) In Example \ref{exaBM}, we have $\aa=k\one$, $S=k[z^\pm
1]$, $m=4$, $\qb=\bar{1}$ and $u=z$. Let $d=1\otimes zd/dz$. We
compute $\varphi(d)(1\otimes z^2)$ from (\ref{inversemap}) with
$\ib=\bar{0}$, $\bar{s}=\bar{2}$, $a_{\ib}=1$, $b_{-\ib+\sb}=z^2$,
$\epi=0$ and $\eps=2$. Then
\begin{eqnarray*} \varphi(d)(1\otimes
z^2)&=& d(1\otimes z^{-2}z^2)z^2+ 2(4)^{-
1}z^{0}[z^{-4}d(1\otimes z^{-0+4})-d(1\otimes z^{0})]z^2\\
&=&0+2(4)^{- 1}[4(1\otimes 1)+0]z^2\\
&=&2(1\otimes z^2).
\end{eqnarray*}
Replacing $2$ with $3$ in the above computations, gives
$\varphi(1\otimes z^3)=3(1\otimes z^2)$. Finally, to compute
$\varphi(1\otimes z^5)$, we note that $z^5\in
S_{\bar{1}}=S_{\0b+\bar{1}}$, so $\sb=\bar{1}$ and $\eps=1$.
Therefore
\begin{eqnarray*} \varphi(d)(1\otimes
z^5)&=& d(1\otimes z^{-1}z^5)z^1+ (4)^{-
1}z^{0}[z^{-4}d(1\otimes z^{-0+4})-d(1\otimes z^{0})]z^5\\
&=&4(1\otimes z^5)+(4)^{- 1}[4(1\otimes 1)+0]z^5\\
&=&5z^5.
\end{eqnarray*}
Thus
$$\varphi(d)(1\otimes z^2z^3)=(1\otimes z^2)\varphi(d)(1\otimes z^3)+
\varphi(d)(1\otimes z^2)(1\otimes z^3).
$$

(ii) Let $\aa=k\one$, $S=k[z^{\pm 1}]$, $\sg_1=\one$ and
$\sg_2(z^n)=\omega^{-n}z^n$. Then
$\big(\dd(S)\big)_{\0b}=\span_{k}\{z^{nm+1}d/dz\mid n\in\bbbz\}$
and $\dd(S_{\0b})=\span_{k}\{t^{n+1}d/dt\mid n\in\bbbz\}$, where
$t=z^m$. Now (\ref{inversemap}) with $u=z^{-1}\in S_{\bar{1}}$
gives the explicit formula for the isomorphism
$\dd(S_{\0b})\cong\big(\dd(S)\big)_{\0b}$. Namely if $d=
t^{n+1}d/dt$, $j\in\bbbz$ and $\eta(j):=m^{-1}(j+\ep(-\jb))$, then
identifying $1\otimes a$ with $a$, we obtain
\begin{eqnarray*}
\varphi(d)(z^j)&=&d(z^{m\eta(j)})z^{-\ep(-\jb)}+\ep(-\jb)m^{-1}[z^{m}d(z^{-m})-d(1)]z^{j}\\
&=& d(t^{\eta(j)})z^{-\ep(-\jb)}+\ep(-\jb)m^{-1}d(t^{-1})z^{m+k}\\
&=&m^{-1}jz^{nm+j}\\
&=&m^{-1}z^{nm+1}d/dz(z^j).
\end{eqnarray*}
So $\varphi$ takes $t^{n+1}d/dt$ to $m^{-1}z^{nm+1}d/dz.$

(iii) Let $\charr(k)=0$ and assume that $\omega$ is a
$m^{\hbox{th}}$-primitive root of unity in $k$. Let $\mathfrak{g}$
be a twisted affine Kac--Moody Lie algebra over $k$. Then
$\mathfrak{g}$ can be realized as the fixed points of a loop
algebra $\aa\otimes S$, under a finite order automorphism
$\sg=\sg_1\otimes\sg_2$, where
\begin{itemize}
\item[-] $\aa$ is a finite dimensional simple Lie algebra,
\item[-] $S=k[z^{\pm 1}]$ is the algebra of Laurent polynomials in
$z$, \item[-] $\sg_1$ is a period $m$ automorphism of $\aa$ and
$\sg_2(z^n)=\omega^{-n}z^n$.
\end{itemize}
Then with respect to the $\zm$-grading on $\aa\otimes S$, we have
$\mathfrak{g}=(\aa\otimes S)_{\bar{0}}$. We also have
$S_{\ib}=z^{-i}k[z^{\pm 1}]$, $C(\aa)=C(\aa)_{\bar{0}}=k$ and as
above $\dd(S)_{\bar{0}}=\span_{k}\{z^{nm+1}d/dz\mid n\in\bbbz\}$.
Therefore, by Theorem \ref{thm2},
\begin{eqnarray*}
\dd(\mathfrak{g})=\dd\big((\aa\otimes
S)_{\0b}\big)&=&\big(\dd(\aa\otimes
S)\big)_{\0b}=\sum_{\ib\in\zm}\dd(\aa)_{\ib}\otimes
S_{-\ib}\soplus
1\otimes\big(\dd(S)\big)_{\0b}\\
&=&\sum_{\ib\in\zm}\ad\aa_{\ib}\otimes z^{-i}k[z^{\pm m}]\soplus
1\otimes\big(\dd(S)\big)_{\0b}.
\end{eqnarray*}
Now it is easy to see that for any $a\in\aa_{\ib}$ and
$n\in\bbbz$, $\varphi(\ad(a)\otimes z^{nm-i})=\ad(a)\otimes
z^{nm-i}$ and $\varphi(1\otimes z^{nm+1}d/dz)=(1\otimes
z^{nm+1}d/dz).$ }
\end{exa}





\end{document}